\documentstyle[12pt]{article}
\begin{document}

\begin{center}
    {\bf \large Gromov-Witten Invariants of Blow-ups Along Points and
                 Curves}
\end{center}
\begin{center} Jianxun Hu \footnote[1]{supported by NNSF of
               China and Lingnan Foundation}

\end{center}

\begin{center}\small Department of Mathematics, Zhongshan University,\\
                Guangzhou, 510275 P. R. China\\
                   and\\
               Department of Mathematics, University of Wisconsin-Madison\\
                 Madison, WI 53706, USA
\end{center}

\begin{center}
\begin{minipage}{130mm}
\vskip 0.5cm
\begin{center}
     {\bf Abstract}
\end{center}
    In this paper, usng the gluing formula of Gromov-Witten invariants under
symplectic cutting, due to Li and Ruan, we studied the Gromov-Witten invariants 
of blow-ups at a smooth point or along a smooth curve. We established some 
relations between Gromov-Witten invariants of $M$ and its blow-ups at a smooth
point or along a smooth curve.

\end{minipage}
\end{center}

\section{Introduction}

 During last several years, there was a great deal of activities to
establish the mathematical foundation of the theory of quantum cohomology
or Gromov-Witten invariants. Ruan-Tian \cite{R1},\cite{RT1},\cite{RT2}
first established for semipositive symplectic manifolds. Recently,
semipositivity condition has been removed by many authors \cite{B},
\cite{FO}, \cite{LT1}, \cite{LT2}, \cite{R2}, \cite{S}. The focus
now is on the calculations and applications. Many Fano manifolds were
computed. We think it is important to study the change of Gromov-Witten 
invariants under surgery. Li-Ruan \cite{LR} gave a gluing formula about 
contact surgery and symplectic cutting. Ionel-Parker \cite{IP} also
studied the Gromov-Witten invariants of symplectic sums. 

Let $\tilde{M}$ be the blow-up of symplectic manifold $M$. There are at 
least two motivations to study the Gromov-Witten invariants of blowups. 
First at all, the curves in the blowup $\tilde{M}$ of a symplectic
manifold $M$ are closely related to curves in $M$. At least for
irreducible curves not contained in the exceptional divisor, we can give a
correspondence between curves in $\tilde{M}$ of a specified homology class 
and curves in $M$ intersecting the blow-up submanifold with a given
multiplicity in terms of the strict transform of curves. Secondly, some
recent research indicated that there is a deep amazing relation between
quantum cohomology and birational geometry. The quantum minimal model 
conjecture, \cite{R3} \cite{R4}, lead to attempt to find quantum
cohomology of a minimal model without knowing minimal model. This problem
requires a thorough understanding of blow-up type formula of Gromov-Witten
invariants and quantum cohomology.

According to McDuff \cite{M1} the blow-up operation in symplectic geometry
amounts to a removal of an open symplectic ball followed by a collapse of
some boundary directions. Lerman \cite{L} gave a generalization of blow-up
construction, `` the symplectic cut". In the case of symplectic manifolds
with hamiltonian circle action, the construction allows us to embedd the
reduced spaces in a symplectic manifold as codimension 2 symplectic
submanifolds.

In this paper, we use symplectic cutting to construct blow-ups at a smooth
point or along smooth submanifolds and use the gluing formula of
Gromov-Witten invariants in \cite{LR} to study the Gromov-Witten
invariants of blow-ups.

Throughout this paper, let $M$ be a compact symplectic manifold of
dimension $2n$, $\tilde{M}$ be the blow-up of $M$ at a smooth point or
along smooth submanifolds. Denote by $p: \tilde{M}\rightarrow M$ the
natural projection. Denote by $\Psi^M_{(A,g)}(\alpha_1,\ldots,\alpha_m)$
the genus $g$ Gromov-Witten invariants of $M$, $\Psi^M_A(\alpha_1, \ldots,
\alpha_m)$ the genus zero Gromov-Witten invariants of $M$.  In this
paper, we established some relations between Gromov-Witten invariants of
$M$ and $\tilde{M}$. Since those curves representing a homology class in
the exceptional divisor have to be contained in the exceptional
divisor and the fact that a GW-invariant
$\Psi^{\tilde{M}}_{(A,g,m)}(\alpha_1,\ldots,\alpha_m) = 0$ if there is
no stable $J$-holomorphic map representing the class $A$ satisfying the
condition given by chomology classes $\alpha_1$, $\ldots$, $\alpha_m$, we
have

{\bf Lemma 1.1:} Suppose that at least one of $\alpha_i$,$ 1\leq i \leq
m$,is the pullback of a cohomology class in $M$ and let $A=re$. Then 
$$
   \Psi^{\tilde{M}}_A(\alpha_1,\ldots,\alpha_m) = 0,
$$
where $e$ denotes the class of a line in the exceptional divisor.

Intuitively, those curves in $\tilde{M}$ which do not intersect with the
exceptional divisor can be identified with curves in $M$. Since GW-invariants 
count curves which represent the given homology class and satisfy the conditions 
given by some cohomology classes, the corresponding invariants on $M$ and 
$\tilde{M}$ should be equal. We showed

{\bf Theorem 1.2:} Suppose that $A\in H_2(M)$, $\alpha_1,\ldots,\alpha_m
\in
H^*(M)$, $g\leq 1$. Then
$$
        \Psi^M_{(A,g)}(\alpha_1.\ldots,\alpha_m) 
     = \Psi^{\tilde{M}}_{(p!(A),g)}(p^*\alpha_1,\ldots,p^*\alpha_m),
$$
where $p!(A) = PDp^*PD(A)$.

We conjecture that the genus condition $g\leq 1$ is a technical one, i.
e. this theorem is true for any genus. But I can not prove it now. I will
study this problem in the future. If we only consider the symplectic 
manifolds of dimension less than $6$, we may remove this condition and
prove

{\bf Theorem 1.3:} Suppose that $\dim_{\bf R}M\leq 6$ and $A\in H_2(M)$, 
$\alpha_1,\ldots , \alpha_m \in H^*(M)$. Then for any genus $g$
$$
       \Psi^M_{(A,g)}(\alpha_1,\ldots,\alpha_m)
      = \Psi^{\tilde{M}}_{(p!(A),g)}(p^*\alpha_1,\ldots,p^*\alpha_m).
$$

From the point of geometry, if we want to express the condition that
curves of homology $A$ pass through a generic point in $M$, we expect to
be able to do this in two different ways: either we add the cohomology
class of a point to the invariants in $M$, or we blow up the point and
count curves with homology class $p!(A)-e$, where $e$ is a class of a line
in the exceptional divisor. We show that these two methods will always
give the same result.

{\bf Theorem 1.4:} Suppose that $A\in H_2(M)$, $\alpha_i \in H^*(M)$,
$1\leq
i \leq m$. Then
$$
        \Psi^M_A (\alpha_1, \ldots ,\alpha_m, [pt])
      = \Psi^{\tilde{M}}_{p!(A)-e}(p^*\alpha_1, \ldots, p^*\alpha_m),
$$
where $e$ denotes the homology class of a line in the exceptional divisor.

So far we are only concerned with the blow-ups of a symplectic manifold at a
smooth point in previous theorems whose proofs are contained in section $3$.
In the case of convex projective variety and genus zero, Gathmann \cite{G}
obtained semilar results to our previous theorems using completely
different method. But the convex condition is very restrictive one and 
most symplectic manifolds are not convex. 

In section $4$, we will discuss the case of blow-ups of any symplectic 
manifolds along some submanifolds. If the blown-up submanifolds are smooth 
curves with nonzero genus or special surfaces, we can show the following 
theorems:

{\bf Theorem 1.5:} Suppose that $C$ is a smooth curve in $M$ such that
either its genus $g_0 \geq 1$ or $g_0 = 0$ and $C_1(M)(C) \geq 0$, 
where $C_1(M)$ denotes the first Chern classes of $M$
and its normal bundle respectively. $A\in H_2(M)$, $\alpha_i\in H^*(M)$,
$1\leq i \leq m$. Then
$$
        \Psi^M_A(\alpha_1, \ldots , \alpha_m)
      = \Psi^{\tilde{M}}_{p!(A)}(p^*\alpha_1, \ldots ,p^*\alpha_m)
$$

About the changes of GW-invariants of blow-up of symplectic manifold
along a smooth surface, in this paper, we only consider the case that
the smooth surface $S$ satisfies one of the followings:
\begin{enumerate}
  \item[(1)] $S = C_1\times C_2$, where $C_1$ and $C_2$ have positive 
genus;
  \item[(2)] $S$ is a $K3$ surface or a torus.
\end{enumerate}

{\bf Theorem 1.6:} If $S$ is a smooth surface in $M$ satisfying one of 
the above two conditions, $A\in H_2(M)$, $\alpha_i\in H^*(M)$, $1\leq 
i \leq m$, satisfy either $deg \alpha_i \geq 2$ or $deg \alpha_i \leq 2$ 
and support away from $S$. Then
$$
         \Psi^M_A (\alpha_1, \ldots, \alpha_m) 
        = \Psi^{\tilde{M}}_{p!(A)}(p^*\alpha_1, \ldots , p^*\alpha_m).
$$

{\bf Acknowledgement:}The author is grateful to the Department of
Mathematics, University of Wisconsin-Madison for its kind hospitality.
The author would like to thank Prof. Yongbin Ruan for his many suggestive 
discussion and encouragement. Thanks also to Prof. An-Min Li, Wanchuan
Zhang, Bohui Chen, Shengda Hu for the valuable dicussions. 

\section{Preliminary Results}

In this section, we describe some notations and preliminary results
that will be used throughout this work. The readers can find their proofs
in the reference \cite{LR}.

\subsection{Symplectic Cutting}

During the last ten years, symplectic surgeries have been sucessfully used
to study symplectic topology, for example, symplectic blow-up and
blow-down by McDuff \cite{MS1} and symplectic norm sum by Gompf \cite{Go}
and McCarthy and Wolfson \cite{MW}. Now we will briefly describe Lerman's 
generalization of the blow-up construction, ``the symplectic cut", \cite{L} 
and \cite{LR}.

Suppose that $ H: M\rightarrow R$ is a periodic hamiltonian function. The
hamiltonian vector field $X_H$ generates a circle action. By adding a
constant, we can assume that $0$ is a regular value. Then, $N = H^{-1}(0)$ 
is a smooth submanifold preserved by circle action. The quotient
$H^{-1}(0)/S^1$ is the famous symplectic reduction. Namely, it has an
induced symplectic structure. Let
$$
         \pi : H^{-1}(0) \longrightarrow Z = H^{-1}(0)/S^1. \eqno(2.1)
$$
$Z$ admits a natural symplectic structure $\tau_0$ such that 
$$
       \pi^*\tau_0 = i^*_0\omega,
$$
where $i_0: H^{-1}(0)\longrightarrow M$ is the inclusion. We note that 
$Z$ is a symplectic orbitfold in general. Furthermore, it is enough that
$H$ is defined in a neighborhood of $H^{-1}(0)$. $(2.1)$ is a circle
bundle.

According to McDuff \cite{M1}, McCarthy-Wolfson \cite{MW}, since $0$ is a
regular value, there is a small interval $I = (-\varepsilon, \varepsilon)$
of regular values. We use a $S^1$-invariant connection on the fibration
$H^{-1}(I)\longrightarrow I$ to show that there is a $S^1$-diffeomorphism
$H^{-1}(I)\cong N\times I$. We will identify $H^{-1}(I)$ with $N\times I$
without any confusion. Then the hamiltonian function is simply the
projection onto the second factor. In such way, we also identify the
symplectic reduction $H^{-1}(t)/S^1$ with $Z$. Suppose that its symplectic
form is $\tau_t$. A beautiful theorem of Duistermaat-Heckman \cite{DH}
says that
$$
    [\tau_t] = [\tau_0] + tc, \eqno(2.2)
$$
where $c$ is the first Chern class of circle bundle $(2.1)$. Hence, if the
boundary components of two symplectic manifolds have the same $\tau_0$,
$c$, we can glue them together.

In the rest of this subsection, we will discuss $\varepsilon$-blow-up
along a submanifold and how to cut the symplectic manifold along a
hypersurface $N$ and collapse the $S^1$-action on $N$ to form two closed
symplectic manifolds if $H^{-1}(I)\cong N\times I$ is symplectically
embedded in a symplectic manifold.

Let $S$ be a compact symplectic submanifold in $(M,\omega)$ of codimension
$2k$. By symplectic neighborhood theorem, there is a tubular neighborhood
${\cal N}_\delta(S)$ of $S$ which is symplectomorphic to the normal bundle
$N_S$. The normal bundle $N_S$ is also a symplectic vector bundle and has
a compatible complex structure. Therefore, we may consider it as a bundle
with fiber $(C^k, -\sqrt{-1}\sum dz_i\wedge d\overline{z_i})$. Furthermore, 
we may consider $N_S$ over $S$ with the symplectic form
$$
       \omega_S = \omega\mid_S + -\sqrt{-1}\sum dz\wedge d\overline{z_i},
$$
where $\omega\mid_S$ is the restriction of the symplectic form $\omega$
to $S$, $z = (z_1, \ldots , z_k)$ are the coordinates in the fiber. The
hamiltonian function is
$$
          H(x, z) = |z|^2 - \varepsilon
$$ 
and the $S^1$-action is given by
$$
         e^{i\theta}(x, z) = (x, e^{i\theta}z).
$$

Consider the symplectic vector bundle $N_S\oplus {\cal O}$ with symplectic
form $\omega_S + -\sqrt{-1}dw\wedge d\bar{w} $ and the momentum map
$\mu (x, z, w) = H(x, z) + |w|^2 $ arising from the actionof $S^1$
on $N_S\oplus {\cal O}$. As  \cite{LR} and  \cite{L}, the manifold 
$M^+ := \{ (x,z)\mid H(x,z) < 0\}$ embeds as an open dense submanifold
into the reduced space
$$
     \overline{M}^+_S := \{(x,z,w)\mid |z|^2 + |w|^2 = \varepsilon\}/S^1
$$
and the difference $\overline{M}^+_S - M^+_S$ is symplectomorphic to the
reduced space $H^{-1}(0)/S^1$.

 A similar procedure defines
$$
    \overline{M}^- := \{(x,z,w)| |z|^2 - |w|^2 = \varepsilon\} /S^1.
$$
It is easy to see that the symplectic manifold $H^{-1}(0)/S^1$ is embedded
on both $\overline{M}_S^+ $ and $\overline{M}_S^-$ as a codimension $2$
symplectic submanifold but with opposite normal bundles. So the symplectic
gluing of $\overline{M}_S^+$ and $\overline{M}_S^-$ along the reduced
space $H^{-1}(0)/S^1$ recovers the neighborhood ${\cal N}_\delta (S)$, i.
e. the normal bundle $N_S$.

We define $\overline{M}^+ := \overline{M}^+_S$ and $\overline{M}^- :=
(M - {\cal N}_\delta (S))\bigcup \overline{M}^-_S$. From the above
description, we know the symplectic gluing of $\overline{M}^+$ and 
$\overline{M}^-$ recovers the original manifold $M$. We will call the
operation that produces $\overline{M}^+$ and $\overline{M}^-$ symplectic
cutting.

Accordign to \cite{MS1},\cite{L},\cite{LR}, we have $\overline{M}^+ =
{\bf P}(N_S\oplus {\cal O})$ and $\overline{M}^- = \tilde{M}$. Specially, 
when $S$ is a point in $M$, we have $\overline{M}^+ = {\bf P}^n$,
$\overline{M}^- = \tilde{M}$.

\subsection{Moduli Spaces}

  Let $(M,\omega)$ be a compact symplectic manifold of dimension $2n$,
$  H : M\longrightarrow R$ a local hamiltonian function such that there is
a small interval $I = (-\delta, \delta)$ of regular values. Denote $N =
H^{-1}(0)$. Suppose that the hamiltonian vector field $X_H$ generates
a circle action on $H^{-1}(I)$. We identify $H^{-1}(I)$ with $I\times N$. 
By a uniqueness theorem on symplectic forms, see \cite{MW}, we may assume 
that the symplectic form on $N\times I$ is expressed by
$$
        \omega = \pi^*(\tau_0 + t \Omega) - \alpha\wedge dt
$$
where $\Omega := d\alpha$ is the curvature  form, which is a $2$-form on
$Z$. We assume that the hypersurface $N = H^{-1}(0)$ divides $M$ into two
parts $M^+$ and $M^-$. As in \cite{LR}, we may consider $M^\pm$ as a
manifold with cylindrical end:
$$
        M^+ = M^+_0 \bigcup \{[0,\infty)\times N\}
$$
$$
        M^- = M^-_0 \bigcup \{(-\infty,0]\times N\}
$$
with symplectic forms $\omega_{\phi^\pm}\mid_{M^\pm_0} = \omega$ and over
the cylinder 
$$
         \omega_{\phi^\pm} = \pi^*(\tau_0 + \phi^\pm\Omega) +
                                 (\phi^\pm)'\alpha \wedge da            
        \eqno (2.3)
$$
where $\phi^+ : [1, \infty) \longrightarrow [-\delta_0, 0)$ and $\phi^- :
(-\infty, -1] \longrightarrow (0,\delta_0]$ are functions such that 
$$
                   \phi^\pm > 0,
$$
$$
        \phi^+(1) = -\delta_0, \,\,\,\, \lim_{a\rightarrow\infty }
                                        \phi^+(a) = 0,
$$             
$$
       \phi^-(-1) = \delta_0, \,\,\,\, \lim_{a\rightarrow
                        -\infty}\phi^-(a) = 0.
$$

For any $J$-holomorphic curve $u : \Sigma\longrightarrow M^\pm$ we 
define the energy $E(u)$ as
$$
     E_{\phi^\pm} = \int_{\Sigma}u^*\omega_{\phi^\pm}.
$$
For any $J$-holomorphic curve $u:\Sigma \longrightarrow {\bf R}\times N$
we write $u = (a,\tilde{u})$ and define
$$
      \tilde{E}_{\phi^\pm}(u) = \int_{\Sigma}\tilde{u}^*(\pi^*\tau_0).
$$
where $\pi$ is the projection in $(2.1)$. Let $(\Sigma, i)$ be a compact
Riemannian surface and $P\subset \Sigma$ be a finite collection of points. Denote
$\stackrel{\circ}{\Sigma} = \Sigma\backslash P$. Let $u :
\stackrel{\circ}{\Sigma}\longrightarrow {\bf R}\times N$ be a
$J$-holomorphic curve, i.e. $u$ satisfies 
$$
          du\circ i = J\circ du.
$$
Following \cite{HWZ1} we impose an enery condition on $u$. Let $\delta_1 <
\delta_2$ be two real numbers and $\Phi$ be the set of all smooth
functions $\phi : {\bf R}\longrightarrow [\delta_1, \delta_2]$ satisfying
$$  \begin{array}{l}
    \phi' > 0\\
   \mbox{$\phi(a) \longrightarrow \delta_2$ as $ a \rightarrow\infty$}\\
   \mbox{$\phi(a) \longrightarrow \delta_1$ as $ a \rightarrow\infty$.}
\end{array}
$$
For any $\phi\in \Phi$ we equip the tube ${\bf R}\times N$ with a
symplectic form $d(\phi \lambda)$. We will call such a $u$ a finite 
energy $J$-holomorphic curve if
$$
      \sup_{\phi\in\Phi}\{\int_{\stackrel{\circ}{\Sigma}}u^*d(\phi\lambda)\} <
\infty.
$$

If we collapse the $S^1$-action on $N= H^{-1}(0)$ we obtain symplectic 
cuts $\overline{M}^+$ and $\overline{M}^-$. The reduced space $Z$ is a
codimension $2$ symplectic submanifold of both $\overline{M}^+$ and 
$\overline{M}^-$. We also can view the symplectic cuts $\overline{M}^+$
and $\overline{M}^-$ as the completions of $M^\pm$ with respect to the
metric $\langle , \rangle_{\omega_{\phi^\pm}}$ see \cite{LR}. We also note
that the almost complex structure on $M^\pm$ is invariant. 

 Let ${\cal M}_{g,m}$ be the moduli space of Riemann surfaces of genus $g$
and with $m$ marked points, and $\overline{\cal M}_{g,m}$ its
Deligne-Mumford compactification. Then $\overline{\cal M}_{g,m}$ consists
of all stable curves of genus $g$ and with $m$ marked points. It is
well-known that $\overline{\cal M}_{g,m}$ is a Kahler orbitfold.

Let $(\Sigma;y_1,\ldots,y_m,p_1,\ldots,p_\nu)\in\overline{\cal
M}_{g,m+\nu}$, and $u : \stackrel{\circ}{\Sigma}\longrightarrow M^\pm$ be 
a finite energy $J$-holomorphic curve. Suppose that $u(z)$ converges to a
$k_i$-periodic orbit $x_{k_i}$as $z$ tends to $p_i$. By using the
removable singularity theorem we get a $J$-holomorphic curve $\bar{u}$
from $\Sigma$ into $ \overline{M}^\pm $. Let $A = [\bar{u}(\Sigma)]$. It
is obvious that 
$$
      E_{\phi^\pm}(u) = \omega(A)
$$
which is independent of $\phi^\pm$. For a map $u$ from $\Sigma$ into ${\bf
R}\times N$, we let $A = [\pi u(\Sigma)]$. Then
$$
     \tilde{E}_{\phi^\pm}(u) = \tau_0(A).
$$
We say $u$ represents the homology class $A$.

{\bf Definition 2.1:} Let $(\stackrel{\circ}{\Sigma}; {\bf y},{\bf p})$ be 
a Riemann surface of genus $g$ with $m$ marked points ${\bf y}$ and $\nu$ 
punctured points ${\bf p}$. A relative stable holomorphic map with $\{k_1,
\ldots, k_\nu\}$-ends from $(\stackrel{\circ}{\Sigma};{\bf y},{\bf p})$ 
into $M^\pm$ is an equivalence class of continuous maps $u$ from
$\stackrel{\circ}{\Sigma}'$ into $(M^\pm)'$, modulo the automorphism group
$stb_u$ and the translations on ${\bf R}\times N$, where 
$\stackrel{\circ}{\Sigma}'$ is obtained by joining chains of ${\bf P^1}$s 
at some double points of $\Sigma$ to separate two components, and then
attaching some trees of ${\bf P^1}$s; $(M^\pm)'$ is obtained by attaching
some ${\bf R}\times N$ to $M^\pm$. We call components of
$\stackrel{\circ}{\Sigma}$ principal components and others bubble
components. Furthermore,
\begin{enumerate}
  \item[(1)] If we attach a tree of ${\bf P^1}$ at a marked point $y_i$ or
a punctured point $p_i$, then $y_i$ or $p_i$ will be replaced by a point
different from intersection points on a component of the tree. Otherwise,
the marked points or punctured points do not change;
  \item[(2)] $\stackrel{\circ}{\Sigma}'$ is a connected curve with normal crossings;
  \item[(3)] Let $m_j$ be the number of special points on $\Sigma_j$ which are
nodal points or marked points or punctured points. Then either
$u|_{\Sigma_j}$ is not a constant or $m_j + 2g_j \geq 3$;
  \item[(4)] The restriction of $u$ to each component is $J$-holomorphic;
  \item[(5)] $u$ converges exponentially to some periodic orbits
$(x_{k_1},\ldots, x_{k_\nu})$ as the variable tends to the punctured
points $(p_1,\ldots,p_\nu)$ repectively; 
  \item[(6)] Let $q$ be a nodal point of $\Sigma'$. Suppose $q$ is the
intersection point of $\Sigma_i$ and $\Sigma_j$. If $q$ is a removable
singular point of $u$, then $u$ is continuous at $q$; If $q$ is a
nonremovable singular point of $u$, then $\Sigma_i$ and $\Sigma_j$ are 
mapped into ${\bf R}\times N$. Furthermore, $u|_{\Sigma_i}$ and
$u|_{\Sigma_j}$ converge exponentially to the same periodic orbit of the
Reeb vector field $X$ on $N$ as the variable tend to the nodal point $q$.
\end{enumerate}

If we drop the condition (4), we simply call $u$ a relative stable map. Let
$\overline{\cal M}_A(M^\pm,g,m,{\bf k})$ be the space of the equivalence
class of relative stable holomorphic curves with ends representing the
homology
class $A$, and  $\overline{\cal B}_A(M^\pm,g,m,{\bf k})$ be the space of
stable maps with ends representing the homology class $A$. 

\subsection{The Fredholm Index}

Denote by $W$ one of $\{ M^+, M^-, {\bf R}\times N\}$. For simplicity, we
consider the case $W = M^\pm$, the situation for ${\bf R}\times N$ is the
same. Let $(\Sigma; y_1,\ldots,y_m,p_1,\ldots,p_\nu) \in {\cal
M}_{g,m+\nu}$, $\stackrel{\circ}{\Sigma} = \Sigma - \{p_1,\ldots,p_\nu\}$.
Let $u: \stackrel{\circ}{\Sigma}\longrightarrow W$ be a finite energy
$J$-holomorphic curve. Suppose that $u(z)$ converges to a $k_i$-periodic
orbit $x_{k_i}$ with $k_i\in {\bf Z}$ as $z$ tends to $p_i$. We consider 
the linearization of $\overline{\partial}$-operator
$$
   D_u = D\overline{\partial}_J(u):C^\infty(\Sigma;u^*TW) 
        \longrightarrow \Omega^{0,1}(u^*TW).
$$

Because the operator $D_u$ is not a Fredholm operator, see \cite{D}, \cite{LR},
To recover Fredholm theory we choose a sufficiently samll weight and define
the weighted Sobolev spaces as follows:

For any section $h\in C^\infty (\Sigma;
u^*TW)$ and section $\eta\in \Omega^{0,1}(u^*TW)$ we define the norms
$$
          \| h \|_{1,p,\alpha} = (\int_\Sigma (|h|^p + |\nabla
           h|^p)d\mu)^{\frac{1}{p}} + (\int_\Sigma e^{2\alpha |s|}(|h|^2
           + |\nabla h|^2)d\mu)^{\frac{1}{2}}  \eqno (2.4)
$$
$$
       \| \eta \|_{p,\alpha} = (\int_\Sigma |\eta|^pd\mu)^{\frac{1}{p}}
          + (\int_\Sigma e^{2\alpha |s|}|\eta|^2d\mu)^{\frac{1}{2}}
         \eqno (2.5)
$$
for $p\ge 2$, where all  norms and covariant derivatives are taken with
respect to the metric $\langle , \rangle$ on $u^*TW$ and the metric on
$\Sigma$. Denote 
$$
       {\cal C}(\Sigma; u^*TW) = \{h\in C^\infty(\Sigma; u^*TW);
         \|h\|_{1,p,\alpha} < \infty\},
$$
$$
      {\cal C}(u^*TW\otimes\wedge^{0,1}) = \{\eta \in \Omega^{0,1}(u^*TW);
            \|\eta\|_{p,\alpha} < \infty\}.
$$
Denote by $W^{1,p,\alpha}(\Sigma; u^*TW)$ and
$L^{p,\alpha}(u^*TW\otimes\wedge^{0,1})$ the completions of ${\cal
C}(\Sigma;u^*TW)$ and ${\cal C}(u^*TW\otimes\wedge^{0,1})$ with repect to
the norms $(2.4)$ and $(2.5)$ respectively. Let $h_{i0} =
(b_{i0},\tilde{h}_{i0}) \in \ker L_{i\infty}$. Put $L_\infty =
(L_{1\infty},\ldots,L_{\nu\infty})$, $h_0 = (h_{10},\ldots, h_{\nu 0})$.
We fix a cutoff function $\rho$:
$$
   \rho (s) = \{\begin{array}{cl}
            1, & \mbox{if $|s| \ge L$,}\\
            0, & \mbox{if $|s| \le \frac{L}{2}$}
          \end{array}
$$
where $L$ is a large positive number. We can consider $h_0$ as a vector
field defined in the Darboux coordinate neighborhood we introduced
previously. We put $\hat{h}_0 = \rho h_0$. Then $\hat{h}_0$ is a section
in $C^\infty(\Sigma;u^*TW)$ supported in the tube $\{(s,t)| \,\, |s|\ge
\frac{L}{2}, t\in S^1\}$. Set
$$
   {\cal W}^{1,p,\alpha} = \{ h + \hat{h}_0 | \,\,\, h\in W^{1,p,\alpha},
            \,\,\,\,  h_0\in \ker L_\infty\}.
$$
The operator 
$$
    D_u : {\cal W}^{1,p,\alpha} \longrightarrow L^{p,\alpha}
$$
is a Fredholm operator so long as $\alpha$ does not lie in the spectrum of
the operator $L_{i\infty}$ for all $i= 1, \ldots, \nu$. We thus have a
Fredholm index ${\bf Ind (D_u,\alpha)}$.

 Let $u = (u^+, u^-) : (\Sigma^+,\Sigma^-) \longrightarrow (M^+,M^-)$ be
$J$-holomorphic curves such that $u^+$ and $u^-$ have $\nu$ ends and they
converge to the same periodic orbits at each end. According to our
convention $\Sigma^\pm$ may not be connected. In this case $Ind(D_{u^\pm},
\alpha)$ denotes the sum of indices of its components. Li and
Ruan \cite{LR} proved the following addition formula for operator $D_u$

{\bf Proposition 2.2: (\cite{LR} Theorem 5.14)} Suppose that $\Sigma = 
\Sigma^+\wedge \Sigma^-$ has genus $g$ and $[u(\Sigma)] = A$, Then 
$$
      Ind (D_{u^+},\alpha) + Ind (D_{u^-},\alpha) 
       = 2(n-1)\nu + 2C_1(A) + 2(3-n)(g-1).  \eqno (2.6)
$$

{\bf Remark 2.3: (\cite{LR} Remark 5.15)} Let $u$ be a $J$-holomorphic curve 
from $(\stackrel{\circ}{\Sigma};y_1,\ldots, y_m,p_1,\ldots,p_\nu)$ into 
$M^\pm$ such that each end converges to a periodic orbit. By using the 
removable singularity theorem we get a $J$-holomorphic curve $\bar{u}$ from 
$\Sigma$ into $\overline{M}^\pm$. Therefore, we have a natural identification 
of finite energy pseudo-holomorphic curves in $M^\pm$ and closed
pseudo-holomorphic curves in the closed symplectic manifolds
$\overline{M}^\pm$. Moreover, the operator $D_u$ is identified with the
operator $D_{\bar{u}}$ in a natural way. Under this identificaton, the
condition that $u$ converges to a $k$-multiple periodic orbit at a
punctured point $p$ is naturally interpreted as $\bar{u}$ being tangent to 
$B$ at $p$ with order $k$. Since $\ker L_\infty$ consists of constant
vectors, we can identify the vector fields in ${\cal W}^{1,p,\alpha}_\pm$ along
$u$ with the vector fields in ${\cal W}^{1,p,\alpha}_\pm$ along $\bar{u}$,
the space $L^{p,\alpha}_\pm$ along $u$ is also identified with the space
$L^{p,\alpha}_\pm$ along $\bar{u}$. Thus we have

{\bf Proposition 2.4: (\cite{LR} Proposition 5.16)}
$$
    Ind (D_u, \alpha) = Ind D_{\bar{u}}.
$$

\subsection{Relative Invariants and Gluing Formula}

From previous subsections, we know that $Z$ is a compact, real codimension two
symplectic submanifold of $\overline{M}^+$ ($\overline{M}^-$ respectively). In 
this section, we will recall the definition of relative GW-invariants for the
pair $(\overline{M}^+,Z)$ and state a gluing formula representing the
GW-invariants of $M$ in terms of the relative GW-invariants of
$(\overline{M}^\pm,Z)$, which are due to Li and Ruan \cite{LR}.

First we recall the definition of virtual neighborhood.

{\bf Definition 2.5: } Let ${\cal M}$ be a compact 
topological space. We call $(U,E,S)$ a virtual neighborhood of $\cal M$ if $U$ 
is a finite dimensional oriented V-manifold (not necessarily compact), $E$ is a 
finite dimensional V-bundle on $U$ and $S$ is a smooth section of $E$ such that
$S^{-1}(0) = {\cal M}$. Suppose that ${\cal M}_{(t)} = \bigcup_t{\cal M}_t\times
\{t\}$ is compact. We call $(U_{(t)},S_{(t)},E_{(t)})$ a virtual neighborhood
cobordism if $U_{(t)}$ is a finite dimensional oriented V-manifold with boundary
and $E_{(t)}$ is a finite dimensional V-bundle and $S_{(t)}$ is a smooth section
such that $S^{-1}_{(t)}(0) = {\cal M}_{(t)}$.

Li and Ruan \cite{LR} proved the following theorem.

{\bf Theorem 2.6: (\cite{LR} Section 7.1)} For $\overline{\cal M}_A(M^+,g,m,
{\bf k})$, there exists a virtual neighborhood $(U,E,S)$.

Using the virtual neighborhood we can define the relative GW-invariants. Recall
that we have two natural maps:
$$
           ev : \overline{\cal B}_A(M^+,g,m,{\bf k}) \longrightarrow (M^+)^m
$$ 
defined by evaluating at marked points and 
$$
           P^+ : \overline{\cal B}_A(M^+,g,m,{\bf k}) \longrightarrow Z^\nu.  
$$
defined by projecting to its periodic orbits. To define the relative
GW-invariants, choose a $r$-form $\Theta$ on $E$ supported in a neighborhood of
the zero section, where $r$ is the dimension of the fiber, such that 
$$
     \int_{E_x} i^*\Theta = 1   
$$
for any $x\in U$, where $i$ is the inclusion map $E_x\longrightarrow E$. We 
call  $\Theta$ a Thom form. Now we can define the relative GW-invariant as
follows:

{\bf Definition 2.7:} Suppose that $\alpha_i\in H^*(M^+,{\bf R})$ and 
$\beta_j \in H^*(S_{k_j},{\bf R})$ represented by differential form.
Define the relative
GW-invariants for $(\overline{M}^+,Z)$ as
$$
   \Psi^{(\overline{M}^+,Z)}_{(A,g,m,{\bf k})}(\alpha_1,\ldots,\alpha_m;
         \beta_1,\ldots,\beta_\nu) = \int_U
       ev^*\Pi^m_{i=1}\alpha_i\wedge \Pi^\nu_{j=1}\beta_j\wedge S^*\Theta.
      \eqno (2.7)
$$
Clearly $\Psi = 0$ if $\sum deg(\alpha_i) + \sum deg(\beta_i) \not= ind$.

Now we want to state a general gluing formula representing GW-invariants
of a closed symplectic manifold in terms of relative GW-invariants of its
symplectic cutting.

In \cite{LR}, Li and Ruan showed that one can glue two
pseudo-holomorphiccurves $(u^+,u^-)$ in $M^+$, $M^-$ with the same end
point to obtain a pseudo-holomorphic curve $u$ in $M$. Suppose that the
homology classes of $u^+$,$u^-$,$u$ are $A^+$,$A^-$,$A$. Denote by 
$\overline{M}^+ \cup_Z \overline{M}^-$ the quotien of $M$ by circle 
action on level set $H^{-1}(0)$. Therefore, we have a projection map
$$
   \pi : M\longrightarrow \overline{M}^+\cup_Z\overline{M}^-. \eqno (2.8)
$$
$\pi$ induces a homomorphism
$$
   \pi_* : H_2(M,{\bf Z}) \longrightarrow H_2(\overline{M}^+\cup_Z\overline{M}^-,
       {\bf Z}).
$$
Using Mayer-Vietoris sequence for $(\overline{M}^+,\overline{M}^-, \overline{M}^+
\cup_Z\overline{M}^-)$, $(u^+, u^-)$ defines a homology class $[u^+ + u^-]\in
H_2(\overline{M}^+\cup_Z\overline{M}^-)$. The existence of glued map $u$ 
implies $[u^+ + u^-] = \pi_*([u])$. If $(u^+, u^-)$ is another representative
and glued to $f'$,
$$
   \pi_*([f']) = [u'^+ + u^-] = [u^+ + u^-] = \pi_*([u]).
$$ 
When $\ker\pi_*\not= 0$, $[u]$, $[u']$ could be different from a vanishing cycle
in $\ker\pi_*$.
   
   Let $[A] = A + \ker\pi_*$. Define 
$$
   \Psi_{([A],\ldots)} = \sum_{B\in [A]}\Psi_{(B,\ldots)}.
$$
By the compactness theorem, the summation of right hand side is finite. To abuse
the notation, we use $[A] = A^+ + A^-$ to represent the set of homology classes 
of glued maps.

The Moduli space $\overline{\cal M}_{[A]}(M,m)$ consists of the components indexed 
by the following data:
\begin{enumerate}
\item[(1)] The combinatorial type of $(\Sigma^\pm,u^\pm) :
\{A_i^\pm,g_i^\pm,m_i^\pm,(k_1^\pm,\ldots,k_\nu^\pm)\}, i=1,\ldots,l^\pm$;
\item[(2)] A map $\rho : \{p^+_1,\ldots,p^+_\nu\}\longrightarrow \{p^-_1,\ldots,
p^-_\nu\}$, where $(p^\pm_1,\ldots,p^\pm_\nu)$ denote the punctured points of
$\Sigma^\pm$.
\end{enumerate}
Using the virtual neighborhood technique as in \cite{R2} and \cite{LR}, we can 
define GW-invariants $\Psi_C$ for each component $C$ and we have
$$
  \Psi_{(M,[A],m)} = \sum \Psi_C.
$$

For the GW-invariants $\Psi_C$, Li and Ruan proved 

{\bf Remark 2.8: (\cite{LR} Remark 7.8)} It is easy to see that 
\begin{enumerate}
\item[(i)] For $C=\{A^+,g^+,m^+\}$, we have
$$
   \Psi_C(\alpha_i^+) =
\Psi^{(\overline{M}^+,Z)}_{(A^+,g^+,m^+)}(\alpha_i^+);
   \eqno(2.9)
$$
\item[(ii)] For $C=\{A^-,g^-,m^-\}$, we have 
$$
   \Psi_C(\alpha_i^-) =
\Psi^{(\overline{M}^-,Z)}_{(A^-,g^-,m^-)}(\alpha_i^-).
  \eqno(2.10)
$$
\end{enumerate}

This remark described the contribution of stable $J$-holomorphic
curves which don't go through the middle to the GW-invariants. Now
we want to state a general gluing formula which describes the contribution 
of stable $J$-holomorphic curves which go through the middle. For simplicit, 
we will only state the gluing formula for the component $C=\{A^+,g^+,m^+,k;$ 
$A^-,g^-,m^-,k\}$ For more general components $C$, Ruan \cite{R4} gave
the steps to write the gluing formula. Choose a homology basis
$\{\beta_b\}$ of $H^*(S_{k_b},{\bf R})$. Let $(\delta_{ab})$ be its
intersection matrix.

{\bf Theorem 2.9: (\cite{LR} Theorem 7.10)} Suppose that $\alpha_i^+|_Z
= \alpha_i^-|_Z$ and hence $\alpha_i^+\cup_Z \alpha_i^- \in
H^*(\overline{M}^+\cup_Z \overline{M}^-, {\bf R})$. Let $\alpha_i = \pi^*
(\alpha_i^+\cup_Z\alpha_i^-)$, where $\pi$ is the map in $(2.8)$.
 
For $C=\{A^+,g^+,m^+,k$; $A^-,g^-,m^-,k\}$. we have the gluing formula
$$
  \Psi_C(\alpha_1,\ldots,\alpha_m) 
=k\sum_{a,b}\sum_{i^+,i^-}\delta^{ab}\Psi^{(\overline{M}^+,Z)}
_{(A^+,g^+,m^+,k)}(\alpha^+_{i^+},\beta_a)\Psi^{(\overline{M}^-,Z)}_{(A^-,g^-,
m^-,k)}(\alpha^-_{i^-}, \beta_b)
  \eqno(2.11)
$$
where $\{i^+,i^-\}$ is an divison of $\{1, \ldots, m\}$.

{\bf Remark 2.10:} For the symplectic blow-up, we have $\ker\pi_* = 0$. Therefore
we have $\Psi_{(A,\ldots)} = \Psi_{([A],\ldots)}$.

\section{Blowup at a smooth point}

In this section, we will only consider the case of blowup at a smooth point. 
We will describe the changes of Gromov-Witten invariants under blowup. 
Actually we will give the proofs of Theorems we state in the introduction.

{\bf Proof of Theorem 1.2:} Let $P_0$ be the blow-up point. We
perform the symplectic cutting for $M$ at $P_0$ as in section 2.1. 
We have
$$
  \overline{M}^+ = {\bf P}^n, \,\,\,\,\, 
  \overline{M}^- = \tilde{M}.
$$

We first consider the contribution of each component to
the GW-invariants. Therefore, we consider the component
$$
  C = \{A^+, g^+, m^+, \{k_1,\ldots, k_\nu\};
        A^-, g^-, m^-, \{k_1,\ldots, k_\nu\}\}.
$$
From Proposition $2.2$, we have
$$
  Ind (D_{u^+},\alpha) + Ind (D_{u^-},\alpha)
    = 2(n-1)\nu + 2C_1(A) + 2(3-n)(g-1).
       \eqno (3.1)
$$

According to our convention, $u^\pm: \Sigma^\pm\longrightarrow M^\pm$
may have many connected components $u_i^\pm: \Sigma_i^\pm\longrightarrow
M^\pm$, $i= 1,\ldots, l^\pm$. Suppose $\Sigma_i^\pm$ has arithemetic genus 
$g_i^\pm$, $g^\pm = \sum g_i^\pm$ with $m_i^\pm$ marked points. Note that 
$\overline{M}^+ = {\bf P}^n$. From Remark $2.3$, it is not difficult to 
see that $\bar{u}_i^+$ can be identified as a stable $J$-holomorphic curve 
$h_i^+$ in ${\bf P}^n$. Then from Proposition $2.4$, we have
\begin{eqnarray*}
  Ind (D_{u^+},\alpha) & = & \sum_{i=1}^{l^\pm} Ind (D_{u_i^+},\alpha)
   = \sum_{i=1}^{l^+} Ind D_{\bar{u}_i^+}\\
  & = & 2\sum C_1([h_i^+]) + 2(3-n)(g^+-l^+) + 2\nu - 2\sum k_i.
\end{eqnarray*}
An intersection multiplicity calculation shows $\sum [h_i^+] = \sum k_i e$,
where $e$ is the homology class of a line in ${\bf P}^n$. Hence 
$\sum C_1([h_i^+]) = (n+1)\sum k_i$. Therefore
$$
  Ind (D_{u^+},\alpha) = 2(3-n)(g^+-l^+) + 2n\sum k_i + 2\nu.
$$
Therefore
$$
   Ind (D_{u^-}, \alpha) = 2C_1(A) + 2(3-n)(g-g^++l^+-1)
            + 2(n-2)\nu - 2n\sum k_i.
$$
Since $\alpha_i \in H^*(M)$, we may assume all $\alpha_i$ suport 
away from the neighborhood ${\cal N}_\delta(P_0)$ (see Section 2) 
of the blowup point $P_0$. So we have $\alpha_i^+ = 0$, $1\leq i \leq m$.
Therefore, if $m^+ > 0$, we have for any $\beta_b\in H^*(Z)$,
$$
   \Psi^{(\overline{M}^+, Z)}_{(A^+,g^+,m^+,\{k_1,\ldots, k_\nu\})}
      (\alpha^+_i, \beta_b) = 0.
$$
This implies $\Psi_C = 0$ except $m^- = m$. Now we assume $m^- = m$,
i. e. $m^+ = 0$. On the other hand, if 
$$
  \sum deg \alpha_i \not= 2C_1(A) + 2(3-n)(g-1) + 2m
$$
where $C_1$ denotes the first Chern class of $M$, by the definition 
of the GW-invariants, we have
$$
  \Psi^M_{(A,g)}(\alpha_1,\ldots,\alpha_m) =
  \Psi^{\tilde{M}}_{(p!(A),g)}(p^*\alpha_1,\ldots, p^*\alpha_m) = 0.
$$
We have proved the assertion of the theorem. Therefore, we also assume
$$
  \sum deg \alpha_i = 2C_1(A) + 2(3-n)(g-1) + 2m. \eqno (3.2)
$$
Since $l^+\geq 1$, $\nu >0$, $g^+\leq g\leq 1$, $k_i > 0$, we have
$$
    2(3-n)(l^+-g^+) - 2\sum k_i - 2\nu < 0.
$$
In fact, if $n\geq 3$, this inequality is obvious. If $n=2$, it follows 
from the inequality $ 2l^+ -2\nu -2\sum k_i < 0 \leq 2g^+$ since $l^+ \leq
\nu$. 

Therefore
\begin{eqnarray*}
    \sum deg (\alpha_i^-) & = & 2C_1(A) + 2(3-n)(g-1) + 2m\\
      & > & 2C_1(A) + 2(3-n)(g-1) + 2m\\
      &   & + 2(3-n)(l^+-g^+) - 2\sum k_i - 2\nu \\
      & \geq & 2C_1(A) + 2(3-n)(g-g^++l^+-1) - 2\sum k_i - 2\nu\\
      &   &             + 2(n-1)(\nu - \sum k_i) + 2m \\
      & = & Ind (D_{u^-}, \alpha) + 2m,
\end{eqnarray*}
since $\nu > 0$, $g\leq 1$, $k_i > 0$. Therefore, by the definition of 
relative GW-invariants, we have for any $\beta_b \in H^*(Z)$,
$$
  \Psi^{(\overline{M}^-,Z)}_{(A^-,g^-,m,\{k_1,\ldots,k_\nu\})}
      (\alpha^-_i, \beta_b) = 0.
$$
Therefore, $\Psi_C = 0$ except $C= \{A^-,g,m\}$.

Now it remains to prove 
$$
   \Psi^{\tilde{M}}_{(p!(A),g)}(p^*\alpha_1,\ldots, p^*\alpha_m)
  = \Psi^{(\overline{M}^-,Z)}_{(A^-,g,m)}(\alpha^-_1,\ldots,
\alpha^-_m).
$$

To prove this, we perform the symplectic cutting for $\tilde{M}$.
Note that the divisor $E$ has normal bundle ${\cal O}(-1)$ in $\tilde{M}$.
We choose the symplectic form
$$
   \tilde{\omega} + dz\wedge d\bar{z}
$$
on ${\cal O}(-1)$. Consider the Hamiltonian function $H(x,z) = |z|^2 - 
\epsilon$ with the $S^1$-action given by 
$$
    e^{2\pi it}(x,z) = (x,e^{2\pi it}z).
$$
We perform the symplectic cutting along the hypersurface $N= H^{-1}(0)$
as in section $2.1$. We have
$$
   \overline{\tilde{M}}^+ = {\bf P}({\cal O}(-1) \oplus {\cal O}),\,\,\,\,\,
  \overline{\tilde{M}}^- \cong \tilde{M}.
$$
Now we use the gluing theorem to prove that the contribution of relative stable
$J$-holomorphic curves in $\tilde{M}$ which touch the exceptional divisor
$E$ to the GW-invariant of $\tilde{M}$ is zero. We consider the component
$$
   C = \{ p!(A)^+, g^+, m^+,\{ k_1,\ldots, k_\nu\};
         p!(A)^-, g^-, m^-, \{k_1,\ldots, k_\nu\}\}.
$$
For the support reasons, we have $\Psi_C = 0$ except 
$$
   C = \{ p!(A)^+, g^+, \{k_1,\ldots, k_\nu\};
         p!(A)^-,g^-,m, \{k_1,\ldots, k_\nu\}\}.
$$
From Proposition $2.2$, we have 
$$
  Ind (D_{u^+},\alpha) + Ind (D_{u^-},\alpha) 
   = 2(n-1)\nu + 2C_1(A) + 2(3-n)(g-1).
$$

As in the first part of our proof, we assume $u^\pm : \Sigma^\pm
\longrightarrow M^\pm$ has $l^\pm$ connected components $u_i^\pm :
\Sigma_i^\pm\longrightarrow M^\pm$, $i=1,\ldots,l^\pm$ and $\Sigma_i^\pm$
has arithemetic genus $g_i^\pm$, $g^\pm = \sum g_i^\pm$ with $m_i^\pm$ 
marked points. From Remark $2.3$, it is not difficult to see that 
$\bar{u}_i^+$ can be identified as a stable $J$-holomorphic curve $h_i^+$ in 
$\overline{\tilde{M}}^+ = {\bf P}({\cal O}(-1) \oplus {\cal O})$. Then
from Proposition $2.4$, we have 
$$
  Ind (D_{u^+}, \alpha) = \sum_{i=1}^{l^+} Ind (D_{u_i^+},\alpha)
      = \sum_{i=1}^{l^+} Ind D_{\bar{u}_i^+}.
$$
To caculate $Ind D_{\bar{u}_i^+}$, we need to extend Mori's cone theory 
to cover stable maps. Mori's cone theory tells us that for any algebraic 
manifold $X$ the set
$$
  NE(X) = \{ \sum_ia_iA_i| a_i\geq 0, A_i \mbox{ is represented by a 
        $J$-holomorphic curve}\}
$$
is a closed cone in $H_2(X,{\bf R})$. We have
       
{\bf Claim:} If $A\in H_2(X,{\bf R})$ is represented by stable 
$J$-holomorphic maps, then $A\in NE(X)$.

In fact, suppose that $A \in H_2(X,{\bf R})$ is represented by a stable
$J$-holomorphic map $f : \Sigma\longrightarrow X$ and $\Sigma$ has
$l$ components $\Sigma_i$. Then $f|_{\Sigma_i} : \Sigma_i
\longrightarrow X$ are $J$-holomorphic curves. Therefore, we have $A = 
\sum [f(\Sigma_i)]$. Hence $A \in NE(X)$. So our claim is true.

Now we want to calculate $Ind D_{\bar{u}_i^+}$. Observe that we obtained 
$\overline{\tilde{M}}^+$ from $M$ by performing the symplectic cutting twice.
We also note that $\overline{\tilde{M}}^+$ is independent of the order of
these two symplectic cuttings. Therefore, if we commute the order of 
these two symplectic cuttings, it is easy to see ${\bf P}({\cal O}(-1)
\oplus {\cal O}) \cong \tilde{\bf P}^n$. By Mori's cone theory, we have
$[h_i^+] = a(L-e) + be$, $a\geq 0$, $b\geq 0$, where $L$ is the class of a
line in $\tilde{\bf P}^n$ with $L\cdot E = 1$ and $e$ is the class of a 
line in the exceptional divisor. Let $H$ be the infinite section in 
${\bf P}({\cal O}(-1)\oplus {\cal O})$. Since $H\cdot [h_i^+] = \sum k_j$,
where summation runs over the ends of $u_i^+$. So we have $a=\sum k_j$.
Since $p!(A)\cdot E = 0$, then $E\cdot [h_i^+] = 2a - b = 0$. Therefore, 
$b = 2\sum k_j $, i. e. $[h_i^+] = \sum k_j (L-e) + 2\sum k_je = \sum k_j L 
+ \sum k_j e$. A simple index caculation shows 
\begin{eqnarray*}
   C_1[h_i^+] & = & ((n+1)H - (n-1)E)\cdot [h_i^+]\\
              & = & ((n+1)H - (n-1)E)\cdot \sum k_j (L + e) = (n+1)\sum k_j.
\end{eqnarray*}
Therefore, we have
\begin{eqnarray*}
  Ind D_{\bar{u}_i^+} & = & 2C_1[h_i^+] + 2(3-n)(g_i^+-1) + 2\nu_i - 2\sum k_j\\
    & = & 2(3-n)(g_i^+-1) + 2\nu_i + 2n\sum k_j
\end{eqnarray*}
where $\nu_i$ is the number of ends in $u_i^+$. Therefore   
\begin{eqnarray*}
   Ind (D_{u^+},\alpha) & = & 2(3-n)(g^+-l^+) + 2\nu + 2n\sum k_i,\\
   Ind (D_{u^-},\alpha) & = & 2C_1(A) + 2(3-n)(g-g^++l^+-1)
         + 2(n-2)\nu - 2n\sum k_i.
\end{eqnarray*}
The same argument as in the first part of the proof shows that 
for any $\beta_b \in H^*(Z)$,
$$
   \Psi^{\overline{\tilde{M}}^-}_{(p!(A)^-,g^-,m,\{k_1,\ldots, k_\nu\})}
       ((p^*\alpha_i)^-,\beta_b) = 0.
$$
Therefore, the contribution of $J$-holomorphic curves to the GW-invariant 
is nonzero only if it doesn't touch the exceptional divisor $E$, i. e. 
$C = \{p!(A)^-,g,m\}$. So from the gluing theorem -- Theorem $2.9$, we have 
$$
   \Psi^{\tilde{M}}_{(p!(A),g,m)}(p^*\alpha_1,\ldots, p^*\alpha_m)
   = \Psi^{(\overline{\tilde{M}}^-,Z)}_{(p!(A)^-,g,m)} ((p^*\alpha_1)^-,
        \ldots, (p^*\alpha_m)^-). \eqno (3.3)
$$
However, $\overline{\tilde{M}}^- = \tilde{M} = \overline{M}^-$. Hence
Theorem $1.2$ follows.

{\bf Proof of Theorem 1.3:} Let $P_0$ be the blow-up point. We perform the 
symplectic cutting for $M$ at $P_0$ as in Section $2.1$. We have
$$
  \overline{M}^+ = {\bf P}^n, \,\,\,\,\, \overline{M}^- = \tilde{M}.
$$

We use the same notations and also first consider the contribution of each 
component to the GW-invariants as in the proof of Theorem $1.2$. Consider the 
component
$$
  C = \{ A^+,g^+,m^+,\{k_1,\ldots,k_\nu\}; A^-,g^-,m^-,\{k_1,\ldots, k_\nu\}\}.
$$
Asimilar calculation as in the proof of Theorem $1.2$ shows
\begin{eqnarray*}
   Ind (D_{u^+}, \alpha) & = & 2(3-n)(g^+-l^+) + 2\nu + 2n\sum k_i\\  
   Ind (D_{u^-}, \alpha) & = & 2C_1(A) + 2(3-n)(g-g^++l^+-1)  
                                + 2(n-2)\nu - 2n\sum k_i.
\end{eqnarray*}
The same argument as in the proof of Theorem $1.2$ shows that the contribution 
of the component $C$ to the GW-invariant of $M$ is nonzero only if $C$ is the 
form
$$
   C = \{A^+,g^+,\{k_1,\ldots,k_\nu\}; A^-,g^-,m,\{k_1,\ldots, k_\nu\}\}.
$$
We also assume 
$$
  \sum deg \alpha_i = 2C_1(A) + 2(3-n)(g-1) + 2m.
$$
Otherwise, the theorem is obvious. So we have 
\begin{eqnarray*}
   \sum deg \alpha_i & > & 2C_1(A) + 2(3-n)(g-g^++l^+-1) - 2\sum k_i - 2\nu + 2m\\
                     & \geq & Ind (D_{u^-}, \alpha ) + 2m.
\end{eqnarray*}
We used the conditions $n\leq 3$, $\nu > 0$, $k_i>0$. Therefore, by the 
definition of relative GW-invariants, we have for any $\beta_b \in H^*(Z)$.
$$
  \Psi^{(\overline{M}^-,Z)}_{(A^-,g^-,m,\{k_1,\ldots, k_\nu\})} 
       (\alpha^-_i,\beta_b) = 0.
$$
Therefore, $\Psi_C = 0$ except $C = \{A^-,g,m\}$. From the gluing theorem, 
we have
$$
  \Psi^M_{(A,g)}(\alpha_1,\ldots,\alpha_m)
        = \Psi^{(\overline{M}^-,Z)}_{(A^-,g,m)}(\alpha_1^-,\ldots, \alpha^-_m).
$$
A similar argument as in the proof of Theorem $1.2$ shows
$$
  \Psi^{(\overline{M}^-,Z)}_{(A^-,g,m)}(\alpha^-_1,\ldots, \alpha^-_m)
  = \Psi^{\tilde{M}}_{(p!(A),g)}(p^*\alpha_1,\ldots, p^*\alpha_m).
$$
we omit this argument.  

{\bf Proof of Theorem 1.4:} We perform symplectic cutting at the point
$P_0$. Then we obtain $\overline{M}^+$, $\overline{M}^-$. Without loss of
generality, we may assume the class $[pt]$ with support in a sufficiently
small neighborhood ${\cal N}_\delta (P_0)$ (see section $2$) of the
blow-up point $P_0$. In fact, we may also assume that $[pt]$ with support
in $M^+$ and $\alpha_i$ with support in $M^-$.

As in the proof of the above theorems, for the reasons of support, the
contribution of the component $C$ to the GW-invariants of $M$ is nonzero 
only if $C$ is the form
$$
  C = \{A^+,1,\{k_1,\ldots, k_\nu\}; A^-,m,\{k_1,\ldots,k_\nu\}\}.
$$
From Proposition $2.2$, we have
$$
  Ind (D_{u^+}, \alpha) + Ind (D_{u^-},\alpha) 
     = 2(n-1)\nu + 2C_1(A) + 2n - 6.
$$

We assume that $u^\pm:\Sigma^\pm\longrightarrow M^\pm$ has $l^\pm$ connected 
components $u_i^\pm : \Sigma_i^\pm\longrightarrow M^\pm$, $i=1,\ldots, l^\pm$. 
From Remark $2.3$, it is not difficult to see that $\bar{u}_i^+$ can be 
identified as a stable $J$-holomorphic curve $h_i^+$ in 
$\overline{M}^+ = {\bf P}^n$. Then from Proposition $2.4$, we have
\begin{eqnarray*}
  Ind (D_{u^+},\alpha) & = & \sum_{i=1}^{l^+} Ind (D_{u_i^+},\alpha)
         = \sum_{i=1}^{l^+} Ind D_{\bar{u}_i^+}\\
    & = & 2\sum_{i=1}^{l^+} C_1[h_i^+] + (2n-6)l^+ + 2\nu - 2\sum k_i.
\end{eqnarray*}
The same calculation as in the proof of Theorem $1.2$ shows $\sum C_1[h_i^+] 
= (n+1)\sum k_i$. Therefore,
$$
   Ind (D_{u^+},\alpha) = (2n-6)l^+ + 2\nu + 2n\sum k_i.
$$
Therefore,
$$
   Ind (D_{u^-}, \alpha) = 2 C_1(A) + (2n-6)(1-l^+) + 2(n-2)\nu - 2n\sum k_i.
$$
We assume, without loss of generality,
$$
  \sum deg \alpha_i + 2n = 2C_1(A) + 2n - 6 + 2m + 2.
$$
Otherwise, for dimension reasons, we have
$$
  \Psi^M_A(\alpha_1,\ldots,\alpha_m,[pt]) 
    = \Psi^{\tilde{M}}_{p!(A)-e}( \alpha_1,\ldots,\alpha_m) = 0.
$$
This proves the assertion of the theorem. Therefore,
$$
  \sum deg(\alpha_i) = 2C_1(A) + 2m - 4.
$$

We claim that the contribution of the component to the GW-invariant of 
$M$ is nonzero only if $l^+ = \nu = k =1$. In fact,  from the connectness of 
stable $J$-holomorphic curves, it is easy to see that $l^+\not= 0$. If $l^+ > 1$,
then $\nu >1$, $\sum k_i > 1$. Therefore, we have
$$
   (2n-6)(1-l^+) + 2(n-2)(\nu - \sum k_i) - 4\sum k_i < -4.
$$
Therefore, we have
\begin{eqnarray*}
  \sum deg \alpha_i & = & 2C_1(A) + 2m - 4\\
    & > & 2C_1(A) + (2n-6)(1-l^+) + 2(n-2)\nu - 2n\sum k_i + 2m\\
    & = & Ind (D_{u^-},\alpha) + 2m.
\end{eqnarray*}
Therefore, for any $\beta_b \in H^*(Z)$,
$$
  \Psi^{(\overline{M}^-,Z)}_{(A^-,m,1)}(\alpha_1,\ldots,\alpha_m;
      \beta_b) = 0.
$$
So the contribution of the component $C$ to the GW-invariants of $M$ is 
nonzero only if
$$
   C = \{ L,1,1; A^-, m, 1\},
$$
where $L$ is the class of a line in ${\bf P}^n$.
From Theorem $2.9$, for the dimension reasons, it follows
$$
   \Psi^M_A(\alpha_1,\ldots, \alpha_m, [pt])
  = \Psi^{(\overline{M}^-,Z)}_{(A^-, m, 1)}(\alpha_1,\ldots,\alpha_m,[Z])
   \Psi^{(\overline{M}^+,Z)}_{(L,1,1)}([pt],[pt]_{Z}),
   \eqno (3.4)
$$
where $[Z]= 1 \in H^0(Z)$, $[pt]_{Z}$ is the fundamental class of
the manifold $Z$, and in the proof of this theorem we will denote
$\alpha_i^-$ and $\alpha_i$ by the same symbol if there is no confusion.

Now we want to prove
$$
  \Psi^{(\overline{M}^+,Z)}_{(L,1,1)} ([pt],[pt]_{Z}) = 1. \eqno (3.5)
$$ 

Before we prove $(3.5)$, we first prove the following claim: For any two
general points in ${\bf P}^n$, we have
$$
  \Psi^{{\bf P}^n}_L([pt],[pt]) = 1.  \eqno (3.6)
$$

Let $J_0$ be the standard complex structure on ${\bf P}^n$. From Lemma
$3.5.1$ in \cite{MS2}, it follows that $D_u$ is surjective for any
$J_0$-holomorphic curve $u : {\bf P}^1\longrightarrow {\bf P}^n$. Hence
we do not need virtual neighborhood to calculate this invariant.
By Theorem $5.3.1$ in \cite{MS2} and the definition of GW-invariant, 
$\Psi^{{\bf P}^n}_e([pt],[pt])$ is exactly the number of lines through
two points (see Example $7.3.1$ in \cite{MS2}). Because two points lie on
a unique line in ${\bf P}^n$, we have
$$
  \Psi^{{\bf P}^n}_L ([pt],[pt]) = 1. \eqno (3.7)
$$
If we choose one of two points in $(3.7)$ to be a general point in the
infinite hyperplane ${\bf P}^{n-1}$, it is not difficult to see from $(3.7)$
$$
  \Psi^{{\bf P}^n}_L([pt],[pt]_{{\bf P}^{n-1}}) = 1, \eqno (3.8)
$$
where $[pt]_{{\bf P}^{n-1}}$ means the point belongs to the infinite
hyperplane ${\bf P}^{n-1}$.

In fact, we may identify $Z$ with ${\bf P}^{n-1}$. Therefore, we may
consider $Z$ as an infinite hyperplane in ${\bf P}^n$. By Remark
$2.3$, we have a natural identification of finite energy 
pseudo-holomorphic curves in $M^+$ and closed pseudo-holomorphic curves in
the closed symplectic manifold $\overline{M}^+ = {\bf P}^n$. The
equality $(3.8)$ tell us that there exists only one unreparameterized
pseudo-holomorphic curve through one point in the infinite hyperplane 
${\bf P}^{n-1}$ and one point outside the infinite hypperplane in 
${\bf P}^n$. Therefore, by the definition of relative GW-invariant and
GW-invariant, we have
$$
 \Psi^{(\overline{M}^+,Z)}_{(e,1,1)}([pt],[pt]_{Z}) 
   = \Psi^{{\bf P}^n}_L([pt],[pt]_{{\bf P}^{n-1}}) = 1.
$$
So we proved $(3.5)$.

To prove our theorem, from $(3.4)$, it suffices to prove
$$
 \Psi^{(\overline{M}^-,Z)}_{(A^-,m,1)} (\alpha_1,\ldots,\alpha_m,1)
  = \Psi^{\tilde{M}}_{p!(A)-e}(\alpha_1,\ldots,\alpha_m). \eqno (3.9)
$$

To prove $(3.9)$, we perform the symplectic cutting for $\tilde{M}$. Note
that the exceptional divisor $E$ has normal bundle ${\cal O}(-1)$ in 
$\tilde{M}$. Therefore, we have 
$$
    \overline{\tilde{M}}^+ = {\bf P}({\cal O}(-1)\oplus {\cal O}), \,\,\,\,\,
   \overline{\tilde{M}}^- \cong \tilde{M}.
$$

Now we consider the contribution of relative stable $J$-holomorphic curves in 
$\tilde{M}$ which touch the exceptional divisor $E$ to the GW-invariants of 
$\tilde{M}$. For the support reason, we only consider the component
$$
  C = \{(p!(A)-e)^+,\{k_1,\ldots,k_\nu\}; (p!(A)-e)^-,m,\{k_1,\ldots,k_\nu\}\}.
$$
From Proposition $2.2$, we have
$$
   Ind (D_{u^+},\alpha) + Ind (D_{u^-},\alpha) = 2C_1(A) + 2(n-1)\nu - 4.
$$

We assume $u^\pm : \Sigma^\pm\longrightarrow M^\pm$ has $l^\pm$ connectec 
components $u_i^\pm : \Sigma_i^\pm \longrightarrow M^\pm$, $i = 1,\ldots, l^\pm$.
From Remark $2.3$, it is not difficult to see that $\bar{u}_i^+$ can be
identified as a stable $J$-holomorphic curve $h_i^+$ in $\overline{\tilde{M}}^+
= {\bf P}({\cal O}(-1)\oplus {\cal O})$. Then from Proposition $2.4$, we have
$$
   Ind (D_{u^+},\alpha) = \sum_{i=1}^{l^+} Ind (D_{u_i^+},\alpha)
        = \sum_{i=1}^{l^+} Ind D_{\bar{u}_i^+}
$$
$$
   = 2\sum C_1[h_i^+] + (2n-6)l^+ + 2\nu - 2\sum k_i.   \eqno (3.10)
$$

Now we want to  calculate $C_1[h_i^+]$. Observe that $\overline{\tilde{M}}^+
= \tilde{\bf P}^n$. By Mori's cone theory, we have $[h_i^+] = a(L-e) + be$,
$a\geq 0$, $b \geq 0$, where $L$ is the class of a line in $\tilde{\bf P}^n$
with $L\cdot E = 1$ and $e$ is the class of a line in the exceptional divisor.
Let $H$ be the infinite section in ${\bf P}({\cal O}(-1)\oplus {\cal O})$.
Since $H\cdot [h_i^+] = \sum k_j$, where summation runs over the ends of
$u_i^+$, so we have $a = \sum k_j$. Since $(p!(A)-e)\cdot E = 1$, then 
$E\cdot [h_i^+] = 2a - b = 1$. Therefore, $b = 2\sum k_j - 1$, i. e. 
$[h_i^+] = \sum k_j L + (\sum k_j -1)e$. Therefore, we have
\begin{eqnarray*}
  C_1[h_i^+] & = & [(n+1)H - (n-1)E]\cdot [\sum k_j L + (\sum k_j - 1)e]\\
  & = & (n+1)\sum k_j - (n-1).
\end{eqnarray*}

Plugging in $(3.10)$, we have
\begin{eqnarray*}
   Ind (D_{u^+},\alpha) & = & 2n\sum k_i - 4l^+ + 2\nu\\
   Ind (D_{u^-},\alpha) & = & 2C_1(A) + 4(l^+-1) + 2(n-2)\nu - 2n\sum k_i.
\end{eqnarray*}

We claim that the contribution of the component to GW-invariant of 
$\tilde{M}$ is nonzero only if $l^+ = \nu = k = 1$.  In fact, we have
\begin{eqnarray*}
   \sum deg \alpha_i & = & 2C_1(A) + 2m -4\\
   & \geq & 2C_1(A) + 4(l^+-1) + 2(n-2)\nu - 2n\sum k_i + 2m\\
   & = & Ind (D_{u^-}, \alpha) + 2m
\end{eqnarray*}
The equality holds if and only if 
$$
  4l^+ + 2(n-2)(\nu - \sum k_i) - 4\sum k_i = 0. \eqno (3.11)
$$
It is easy to see $(3.11)$ holds if and only if $l^+ = \nu = \sum k_i$
because $ l^+ \leq \nu \leq \sum k_i$. From $\nu = \sum k_i$ it follws that
$k_i = 1$. Hence Each componet $[h_i^+]$ is just the line $L$. If $l^+
>1$, then we have $ 1 = E\cdot (p!(A) - e) = E\cdot \sum [h_i^+] = l^+$. 
This is a contradiction. So the contribution of the component $C$ to 
the GW-invariants of $\tilde{M}$ is nonzero only if 
$$
  C = \{ L, 1, 1; (p!(A) - e)^-, m,1\}.
$$
From Theorem $2.9$, for the dimension reasons, it follows
$$
  \Psi^{\tilde{M}}_{p!(A) - e}(\alpha_1,\ldots, \alpha_m)
  = \Psi^{(\overline{\tilde{M}},Z)}_{((p!(A)-e)^-,m,1)}(\alpha_1,\ldots,
   \alpha_m, 1)\Psi^{(\overline{\tilde{M}},Z)}_{(L,1)}([pt]_Z). \eqno (3.12)
$$
Because there is a unique line passing through a point in the infinite
section of $\tilde{\bf P}^n$ and intersecting at one point with the exceptional 
dvisor, it is easy to show that 
$$
   \Psi^{(\overline{\tilde{M}}^-,Z)}_{(L,1)}([pt]_Z) = 1.
$$
Therefore, we have
$$
   \Psi^{\tilde{M}}_{p!(A)-e}(\alpha_1,\ldots,\alpha_m)
  = \Psi^{(\overline{\tilde{M}}^-,Z)}_{((p!(A)-e)^-,m,1)}(\alpha_1,
    \ldots, \alpha_m,1). \eqno (3.13)
$$ 

From $(3.9)$ and $(3.13)$, to prove our theorem, it suffices to prove
$$
  \Psi^{(\overline{M}^-,Z)}_{(A^-,m,1)}(\alpha_1,\ldots,\alpha_m,1)
  = \Psi^{(\overline{\tilde{M}}^-,Z)}_{((p!(A)-e)^-,m,1)}(\alpha_1,
   \ldots,\alpha_m,1).  \eqno (3.14)
$$

Assume that $\tilde{u}: \Sigma\longrightarrow \tilde{M}$ is a 
pseudo-holomorphic curve representing $p!(A)-e$. Performing symplectic
cutting, we obtained $\tilde{u}^\pm: \Sigma^\pm\longrightarrow
\overline{\tilde{M}}^\pm$ and $[\tilde{u}^\pm] = (p!(A)-e)^\pm$. Let 
$p : \tilde{M}\longrightarrow M$ be the projection of the blowup. The
map $p\tilde{u} : \Sigma\longrightarrow M$ is also a pseudo-holomorphic
curve representing $A$. Since $\overline{\tilde{M}}^-\cong \tilde{M}
\cong \overline{M}^-$, we may consider $(p!(A)-e)^-$ and $A^-$ as 
homology classes in a same manifold $\tilde{M}$. From the calculation 
in our proof, it follows that $[p\tilde{u}^+]$ is the class of a line
in $\overline{M}^+ = {\bf P}^n$. From Remark $2.10$, we have $[p\tilde{u}^+
+ p\tilde{u}^-] = A$ and $[p\tilde{u}^-] = A^-$. From symplectic cutting,
we may identify $\tilde{u}^-$ and $p\tilde{u}^-$ in $\tilde{M}\setminus E$.
Therefore, $A^- = (p!(A) - e)^-$. By the definition of relative GW-invariant,
we have $(3.14)$. This proves Theorem $1.4$.   
    
{\bf Corollary:} $\Psi^{\tilde{M}}_e([pt]_E,[pt]_E) = 1$, where 
$[pt]_E$ denotes the fundamental class of the exceptional divisor $E$ and
$e$ is the class of a line in the exceptional dvisor $E$.

{\bf Proof:} Lemma $1.1$ tells us that those curves representing a homology 
class in the exceptional divisor have to be contained in the exceptional 
divisor $E$. Since $E$ may be identified with ${\bf P}^{n-1}$. So the 
corollary follows from $(3.6)$.

\section{Blow-up along submanifolds}

In last section, we described some changes of GW-invariants under 
blow-up of symplectic manifold at a general point. In this section, we
will consider the changes of GW-invariants of blow-up of symplectic
manifold along a smooth curve or an smooth surface. As the author knew, so
far only Gathmann \cite{G} delt with two easy examples: the blow-up of a
space curve $Y\subset {\bf P}^3$ and the blow-up of an abelian surface in
${\bf P}^4$.

{\bf Proof of Theorem 1.5:} Since $C$ is a smooth curve of $M$, the normal
bundle $N_C$ is a symplectic vector bundle. By symplectic neighborhood
theorem, there is a tubular neighborhood ${\cal N}_\delta(C)$ of $C$ which
is symplectomorphic to the normal bundle $N_C$. We perform the symplectic
cutting as in section $2.1$. We obtained
$$
  \overline{M}^+ = {\bf P}(N_C\oplus {\cal O}), \,\,\,\,\,
  \overline{M}^- = \tilde{M}.
$$

From the divisor property and skew symmetry of GW-invariants, without loss 
of generality, we may assume that $deg \alpha_i > 2$, $1\leq i\leq m$. 
Therefore, if we choose a sufficiently small $\delta > 0$, we may also assume 
$\alpha^+_i = 0$.

Similar to the proof of Theorem $1.2$, we first consider the contribution of
each component to the GW-invariants. Therefore, we consider the component 
$$
   C = \{A^+,m^+,\{k_1,\ldots,k_\nu\}; A^-,m^-,\{k_1,\ldots,k_\nu\}\}.
$$
From Proposition $2.2$, we have
$$
  Ind (D_{u^+},\alpha) + Ind (D_{u^-},\alpha) = 2(n-1)\nu + 2C_1(A) + 2n - 6.
$$

As in the proof of Theorem $1.2$, we assume $u^\pm : \Sigma^\pm\longrightarrow
M^\pm$ has $l^\pm$ connected components $u_i^\pm : \Sigma_i^\pm\longrightarrow
M^\pm$, $i=1,\ldots,l^\pm$. From Remark $2.3$, it is not difficult to see that 
$\bar{u}_i^+$ can be identified as a stable $J$-holomorphic curve $h_i^+$ in 
$\overline{M}^+ = {\bf P}(N_C\oplus {\cal O})$. Then from Proposition $2.4$, 
we have 
$$
  Ind (D_{u^+},\alpha)  =  \sum_{i=1}^{l^+} Ind (D_{u_i^+},\alpha)
            = \sum_{i=1}^{l^+} Ind D_{\bar{u}_i^+}
$$
$$
    =  2\sum_{i=1}^{l^+} C_1[h_i^+] + (2n-6)l^+ + 2\nu - 2\sum k_i.
    \eqno (4.1)
$$

Now we want to calculate $C_1[h_i^+]$ in two cases of our theorem.

{\bf Case $1$:} The genus $g_0\geq 1$.
 
In this case, we claim that all stable $J$-holomorphic maps $h_i^+$
can only stay in fibers of $\overline{M}^+ = {\bf P}(N_C\oplus {\cal O})$.
Otherwise, suppose that there is a stable $J$-holomrphic curve $h_i^+:
\Sigma\longrightarrow \overline{M}^+$ which doesn't stay in a fiber.
Since we only consider the genus zero GW-invariants, we assume that 
$\Sigma$ has genus zero. Denote by $\pi : {\bf P}(N_C\oplus {\cal O})
\longrightarrow C$ the projection of the projective bundle. Then we 
have a stable $J$-holomorphic map $\pi\circ h_i^+ : \Sigma \longrightarrow
C$ satisfying $[\pi\circ h_i^+]\not= 0$. We can perform pre-gluing as in 
the section $6$ of \cite{LR} and obtain a system of small perturbed 
$J$-holomorphic curves $f_n : \Sigma_n\longrightarrow C$ which represent 
the class $[\pi\circ h_i^+]$ and satisfy the perturbed Cauchy-Riemann 
equation $\overline{\partial}_Jf_n = \nu_n$, here $\Sigma_n$ is a smooth 
Riemann surface. Actually we can choose $\nu_n\longrightarrow 0$ as 
$n\longrightarrow \infty$. Therefore, by Gromov compactness theorem, we 
have that $f_n$ weakly converges to a (possibly reducible) $J$-holomorphic 
curve $u = (u^1,\ldots,u^N)$ and $[\pi\circ h_i^+] = \sum_{j=1}^{N} [u^j]
\not= 0$. Therefore we have a nonconstant $J$-holomorphic curve $f:\Sigma_1
\longrightarrow C$ and $\Sigma_1$ has genus zero. it is wellknown that if 
$f' : S\longrightarrow S'$ is a holomorphic map between compact Riemann 
surfaces, then the genus of $S$ and $S'$ satisfy $g(S)\geq g(S')$ unless 
$f'$ is constant (see \cite{GH} p.$219$).  Since $g(C) = g_0 \geq 1$, we 
have a contradiction. So our claim is true.  

An simple index calculation shows $C_1[h_i^+] = n\sum k_j$
where summation runs over ends of component $u_i^+$. In this case, 
we have
$$
   Ind (D_{u^+},\alpha) = (2n-6)l^+ + 2(n-1)\sum k_i + 2\nu.
$$

{\bf Case $2$:} $g_0 = 0$ and $C_1(M)(C)\geq 0$.

A simple calculation show that $C_1({\bf P}(N_C\oplus {\cal O}))
= C_1(C) + C_1(N_C) + n\xi$ $= C_1(M) + n\xi$, here $\xi $ is the 
class of infinite section in ${\bf P}(N_C\oplus {\cal O})$ over $C$. 
Therefore, from the assumption of the theorem and an intersection 
multiplicity calculation shows 
$$
   \sum_{i=1}^{l^+} C_1[h_i^+] \geq n\sum k_i.
$$
In this case, we have
$$
  Ind (D_{u^+},\alpha) \geq (2n-6)l^+ + 2(n-1)\sum k_i + 2\nu.
$$

Summarise the above two cases, we have
\begin{eqnarray*}
   Ind (D_{u^+},\alpha) & \geq & (2n-6)l^+ + 2(n-1)\sum k_i + 2\nu,\\
   Ind (D_{u^-},\alpha) & \leq & 2C_1(A) + (2n-6)(1-l^+) 
             - 2(n-1)(\nu - \sum k_i) - 2\nu.
\end{eqnarray*}

Since $\alpha^+_i=0$, $1\leq i\leq m$, if $m^+ > 0$, we have for any 
$ \beta_b \in H^*(Z)$
$$
  \Psi^{(\overline{M}^-,Z)}_{(A^+,m^+,\{k_1,\ldots,k_\nu\})}
    (\alpha^+_i, \beta_b) = 0.
$$
This implies $\Psi_C = 0$ except $m^- = m$. So we may assume $m^- = m$.
By the same argument in the proof of Theorem $1.2$, we also may assume
$$
  \sum deg \alpha_i = 2C_1(A) + 2n - 6 + 2m.
$$
Then
\begin{eqnarray*}
  \sum deg (\alpha_i^-) & = & 2C_1(A) + 2n - 6 + 2m\\
         & > & 2C_1(A) + (2n-6)(1-l^+) + 2(n-1)(\nu - \sum k_i) - 2\nu +2m\\
         & \geq & Ind (D_{u^-},\alpha) + 2m^-,
\end{eqnarray*}
since $\nu > 0$, $k_i > 0$, $n\geq 3$. Therefore, by the definition of
relative GW-invariants, we have for any $\beta_b\in H^*(Z)$
$$
  \Psi^{(\overline{M}^-,Z)}_{(A^-,m,\{k_1,\ldots,k_\nu\})}
  (\alpha^-_i,\beta_b) = 0.
$$
Therefore, $\Psi_C = 0$ except $C = \{A^-,g,m\}$.

So now it remains to show 
$$
   \Psi^{\tilde{M}}_{p!(A)}(p^*\alpha_1,\ldots,p^*\alpha_m)
  = \Psi^{(\overline{M}^-,Z)}_{(A^-,m)}(\alpha^-_1,\ldots,\alpha^-_m).
$$
To prove this, we perform the symplectic cutting for $\tilde{M}$ around $E$
as in the proof Theorem $1.2$. Therefore, we have
$$
  \overline{\tilde{M}}^+ = {\bf P}(N_E\oplus {\cal O}),\,\,\,\,\,
  \overline{\tilde{M}}^- \cong \tilde{M}.
$$
Now we use the gluing theorem to prove the contribution of stable 
$J$-holomorphic curves in $\tilde{M}$ which touch the exceptional divisor 
$E$ to the GW-invariants of $\tilde{M}$ is zero. We consider the component
$$
 C =\{p!(A)^+,m^+,\{k_1,\ldots,k_\nu\};p!(A)^-,m^-,\{k_1,\ldots,k_\nu\}\}.
$$
Since $\alpha^+_i = 0$, $1\leq i\leq m$, we have $\Psi_C = 0$ except
$$
  C =\{p!(A)^+,\{k_1,\ldots,k_\nu\};p!(A)^-,m,\{k_1,\ldots,k_\nu\}\}.
$$
From Proposition $2.2$, we have 
$$
  Ind (D_{u^+},\alpha) + Ind (D_{u^-},\alpha) = 2(n-1)\nu + 2C_1(A) + 2n-6,
$$
where $C_1$ denotes the first Chern class of $M$.

We assume that $u^\pm : \Sigma^\pm\longrightarrow M^\pm$ has $l^\pm$
connected components $u_i^\pm : \Sigma_i^\pm\longrightarrow M^\pm$,
$ i=1,\ldots, l^\pm$. From Remark $2.3$, it is not difficult to see
that $\bar{u}_i^+$ can be identified as stable $J$-holomorphic curve 
$ h_i^+$ in $\overline{\tilde{M}}$. then from Proposition $2.4$, 
we have
$$
  Ind (D_{u^+}, \alpha)  =  \sum_{i=1}^{l^+} Ind (D_{u_i^+},\alpha)
      = \sum_{i=1}^{l^+} Ind D_{\bar{u}_i^+}
$$
$$  
      =  (2n-6)l^+ + 2\sum_{i=1}^{l^+}C_1[h_i^+] + 2\nu - 2\sum k_i,
   \eqno (4.2)
$$
where $C_1$ is thefirst Chern class of $\overline{\tilde{M}}^+$.

Let $V$ be a complex rank $r$ vector bundle over $X$, and $\pi :{\bf P}(V)
\longrightarrow X$ be the corresponding projective bundle. Let $\xi_V$ be 
the first Chern class of the tautological line bundle in ${\bf P}(V)$. A 
simple calculation shows 
$$
  C_1({\bf P}(V)) = \pi^* C_1(X) + \pi^*C_1(V) - r\xi_V.  \eqno (4.3)
$$

Note that $\overline{\tilde{M}}^+ = {\bf P}(N_E\oplus 
{\cal O})$ and $ E = {\bf P}(N_C)$. Applying $(4.3)$ to 
$\overline{\tilde{M}}^+$ and $E$, we obtain
\begin{eqnarray*}
  C_1(\overline{\tilde{M}}^+) & = & C_1(E) + C_1(N_E) - 2\xi\\
   & = & C_1(C) + C_1(N_C) - (n-1)\xi_1 + C_1(N_E) + 2\xi,
\end{eqnarray*}
where $\xi_1$ and $\xi$ are the first Chern classes of the tautological
line bundles in ${\bf P}(N_C)$ and ${\bf P}(N_E\oplus {\cal O})$
respctively. Here we denote Chern class and its pullback by a same symbol.
It is wellknow that the normal bundle to $E$ in $\tilde{M}$ is just the 
tautological bundle on $E\cong {\bf P}(N_C)$. Therefore $C_1(N_E) = \xi_1$. 
So we have
$$
  C_1(\overline{\tilde{M}}^+) = C_1(M) - (n-2)\xi_1 - 2\xi.
$$

We know that $\overline{\tilde{M}}$ is a projective bundle over $E$ 
with fiber ${\bf P}^1$. Let $L$ be the class of a line in the fiber 
${\bf P}^1$ and $e$ be the class of a line in the fiber ${\bf P}^{n-2}$
in $E={\bf P}(N_C)$. Denote by $[h_i^+]^C$ the homology class of the 
projection in $C$ of the curve $h_i^+$. Denote by $[h_i^+]^F$ the 
difference of $[h_i^+]$ and $[h_i^+]^F$ i. e. $[h_i^+]^F = [h_i^+] -
[h_i^+]^C$. Then it is easy to know $[h_i^+]^F = aL + be$. Since 
$\xi\cdot [h_i^+] = \sum k_j $, where the summation runs over ends of 
$u_i^+$, and $E\cdot [h_i^+] = 0$, so we have $\xi\cdot [h_i^+]^F = a = 
\sum k_j$ and $E\cdot [h_i^+]^F = a - b = 0$. Therefore, we have $a = b 
= \sum k_j$. So we have $[h_i^+]^F = \sum k_j (L+e)$. For Case 1, we 
have $[h_i^+]^C = 0$. Therefore, we have
$$
   \sum_{i=1}^{l^+} C_1[h_i^+] = 2(n-1) \sum k_i.
$$
For Case $2$, since $C_1(C) + C_1(N_C) \geq 0$, we have 
$$
   \sum_{i=1}^{l^+} C_1[h_i^+] \geq  2(n-1)\sum k_i.
$$
Plugging in $(4.2)$, we have
$$
  Ind (D_{u^+},\alpha) \geq (2n-6)l^+ + 2(2n-3)\sum k_i + 2\nu.
$$
Therefore,
$$
  Ind (D_{u^-},\alpha) \leq 2C_1(A) + (2n-6)(1-l^+) 
        + (2n-2)(\nu - \sum k_i) - 2(n-2)\sum k_i.
$$   

For the same reasons as in the proof of Theorem $1.2$, we also may assume
$$
   \sum deg (p^*\alpha_i) = 2C_1(A) + 2n-6 + 2m.
$$
Then,
\begin{eqnarray*}
  \sum deg (p^*\alpha_i) & = & 2C_1(A) + 2n-6 + 2m\\
      & > & 2C_1(A) + (2n-6)(1-l^+) + (2n-2)(\nu - \sum k_i) \\
      &   &      - 2(n-2)\sum k_i + 2m\\
      & \geq & Ind (D_{u^-},\alpha) + 2m,
\end{eqnarray*}
since $\nu > 0$, $k_i > 0$. Therefore, by the definition of relative 
GW-invariants, we have for any $\beta_b\in H^*(Z)$,
$$
  \Psi^{(\overline{\tilde{M}}^-,Z)}_{(p!(A)^-,m,\{k_1,\ldots,k_\nu\})}
   ((p^*\alpha_i)^-,\beta_b) = 0.
$$
Therefore, the contribution of $J$-holomorphic curves to the GW-invariant
is nonzero only if it doesn't touch the exceptional divisor $E$, i. e. 
$C=\{p!(A)^-,m\}$. So from Theorem $2.9$, we have
$$
  \Psi^{\tilde{M}}_{(p!(A),m)}(p^*\alpha_1,\ldots,p^*\alpha_m)
 = \Psi^{(\overline{\tilde{M}}^-,Z)}_{(p!(A)^-,m)}
   ((p^*\alpha_1)^-,\ldots, (p^*\alpha_m)^-).
$$
The rest of the proof is the same as that of the proof of Theorem $1.2$. So we 
omit it here. This proves Theorem $1.5$.

{\bf Proof of Theorem 1.6:} Since $S$ is a smooth surface, the normal bundle
$N_S$ is a symplectic vector bundle. By symplectic neighborhood theorem, 
there is a tubular neighborhood ${\cal N}_\delta (S)$ of $S$ which is 
symplectomorphic to the normal bundle $N_S$. We perform the symplectic cutting 
as in section $2.1$. We obtain
$$
  \overline{M}^+ = {\bf P}(N_S\oplus {\cal O}), \,\,\,\,\,
   \overline{M}^- = \tilde{M}.
$$
We may assume $\alpha_i^+ = 0$ if we choose a sufficiently small $\delta > 0$
because of the assumption of $\alpha_i$.

Similar to the proof of Theorem $1.5$, we first consider the contribution of 
each component to the GW-invariants. Therefore, we consider the component
$$
  C = \{ A^+,m^+,\{k_1,\ldots,k_\nu\}; A^-,m^-,\{k_1,\ldots,k_\nu\}\}.
$$
From Proposition $2.2$, we have
$$
   Ind (D_{u^+},\alpha) + Ind (D_{u^-},\alpha) = 2(n-1)\nu + 2C_1(A) + 2n-6.
$$

We assume $u^\pm : \Sigma^\pm\longrightarrow M^\pm$ has $l^\pm$ connected 
components $u_i^\pm : \Sigma_i^\pm\longrightarrow M^\pm$, $i=1,\ldots,l^\pm$.
From Remark $2.3$, it is not difficult to see $u_i^+$ can be identified as
a stable $J$-holomorphic curve $h_i^+$ in $\overline{M}^+ = {\bf P}(N_S\oplus 
{\cal O})$. Then from Proposition $2.4$, we have
\begin{eqnarray*}
    Ind (D_{u^+},\alpha) & = & \sum_{i=1}^{l^+} Ind (D_{u_i^+},\alpha)
             = \sum_{i=1}^{l^+} Ind D_{\overline{u}_i^+}\\
   & = & 2\sum_{i=1}^{l^+} C_1[h_i^+] + (2n-6)l^+ + 2\nu - 2\sum k_i,
\end{eqnarray*}
where $C_1$ is the first Chern class of $\overline{M}^+$.

Now we want to calculate $C_1[h_i^+]$. It is wellknown that there is no 
nonconstant stable $J$-holomorphic curves in $S$ if $S$ satisfies the 
conditions (2). If $S$ satisfies the condition (1), the similar argument as 
in the case 1 of Theorem $1.5$ shows there is no nonconstant stable 
$J$-holomorphic curves in $S$. Therefore, all stable $J$-holomorphic curves 
$h_i^+$ can only stay in fibers of $\overline{M}^+ = {\bf P}(N_S\oplus {\cal O})$
over $S$.

An simple index calculation shows $C_1[h_i^+] = (n-1)\sum k_j$ where summation 
runs over ends of component $u_i^+$. In this case, we have
$$
   Ind (D_{u^+},\alpha) = (2n-6)l^+ + 2(n-2)\sum k_i + 2\nu.
$$
Therefore, we have
\begin{eqnarray*}
   Ind (D_{u^-},\alpha) & = & 2C_1(A) + (2n-6)(1-l^+) 
                                + 2(n-2)(\nu - \sum k_i).
\end{eqnarray*}

The same argument as in the proof of Theorem $1.5$ shows that the 
contribution of the component $C$ to the GW-invariant of $M$ is nonzero
only if $S$ is the form
$$
  C = \{ A^+,\{k_1,\ldots, k_\nu\}; A^-,m,\{k_1,\ldots, k_\nu\}\}.
$$
We also assume
$$
   \sum deg \alpha_i = 2C_1(A) + 2n-6 + 2m.
$$
The same argument as in the proof of theorem $1.5$ shows $\Psi_C = 0$
except $C = \{ A^-,m\}$. From the gluing theorem, we have
$$
      \Psi^M_A(\alpha_1,\ldots, \alpha_m)
     = \Psi^{(\overline{M}^-,Z)}_{(A^-,m)}(\alpha_1^-,\ldots, \alpha_m^-).
$$

Now it remains to prove
$$
    \Psi^{\tilde{M}}_{p!(A)}(p^*\alpha_1,\ldots, p^*\alpha_m)
    = \Psi^{(\overline{M}^-,Z)}_{(A^-,m)}(\alpha_1^-,\ldots, \alpha_m^-).
$$

To prove this, we perform the symplectic cutting for $\tilde{M}$ around
$E$ as in the proof of Theorem $1.5$. Therefore, we have
$$
   \overline{\tilde{M}}^+ = {\bf P}(N_E\oplus {\cal O}), \,\,\,\,\,
   \overline{\tilde{M}}^- \cong \tilde{M}.
$$

We also use the gluing theorem to prove that the contribution of stable 
$J$-holomorphic curves in $\tilde{M}$ which touch the exceptional divisor 
$E$ to the GW-invariant of $\tilde{M}$ is zero. We consider the component
$$
  C = \{p!(A)^+,m^+,\{k_1,\ldots,k_\nu\}; p!(A)^-,m^-,\{k_1,\ldots,K_\nu\}\}.
$$

Since $\alpha_i^+ = 0$, $1\leq i\leq m$, we have $\Psi_C = 0$ except
$$
  C = \{p!(A)^+, \{k_1,\ldots,k_\nu\}; p!(A)^-,m,\{k_1,\ldots,k_\nu\}\}.
$$
The similar calculation to that in the proof of Theorem $1.5$ shows
\begin{eqnarray*}
    Ind (D_{u^+},\alpha) & = & (2n-6)l^+ + 2\nu + 2(2n-5)\sum k_i,\\
    Ind (D_{u^-},\alpha) & = & 2C_1(A) + (2n-6)(1-l^+)
                        + 2(n-2)(\nu -\sum k_i) - 2(n-4)\sum k_i.
\end{eqnarray*}
The rest of the proof is the same as that of the proof Theorem $1.5$.
so we omit it. This completes the proof of Theorem $1.6$.

Email address: stsjxhu@zsu.edu.cn

\end{document}